\documentclass[preprint,1p,times,dvipsnames,authoryear]{elsarticle}

\usepackage{amssymb}
\usepackage{amsmath}

\usepackage[usenames,dvipsnames]{color}	
\usepackage[table, dvipsnames]{xcolor}

\newcommand{\om}{K} 
\newcommand{\el}{K}
\newcommand{\domain}{K}

\newcommand{\norm}[1]{\left\lVert#1\right\rVert}

\newcommand{\Cp}{C_p}
\newcommand{\Ci}{C_i}
\newcommand{\Ct}{C_t}

\newcommand{\ie}{i.\,e.\ }

\newcommand{\eg}{e.\,g.,\ }

\newcommand{\CC}{\mathrm{CC}}
\newcommand{\IC}{\mathrm{IC}}
\newcommand{\CR}{\mathrm{CR}}
\newcommand{\AR}{\mathrm{AR}}
\newcommand{\KE}{\mathrm{KE}}
\newcommand{\KAR}{\mathrm{KAR}}
\newcommand{\APR}{\mathrm{APR}}
\newcommand{\SE}{\mathrm{SE}}
\newcommand{\sSE}{\mathrm{sSE}}
\newcommand{\ER}{\mathrm{ER}}
\newcommand{\MPD}{\mathrm{MPD}}
\newcommand{\sMPD}{\mathrm{sMPD}}
\newcommand{\MA}{\mathrm{MA}}
\newcommand{\MX}{\mathrm{MX}}
\newcommand{\SR}{\mathrm{SR}}
\newcommand{\ISO}{\mathrm{ISO}}

\newcommand{\grad}{\nabla}

\definecolor{dgreen}{rgb}{0.0, 0.5, 0.0}

\definecolor{cnrblue}{rgb}{0.21176470588235294, 0.40784313725490196, 0.6392156862745098}

\usepackage{tikz}
\usepackage{pgfplots}
\usepgfplotslibrary{statistics} 

\usepackage{pgfplotstable}
\usepackage{booktabs}
\usepackage{tikz}

\usetikzlibrary{pgfplots.statistics}
\usepackage{tikz}
\usetikzlibrary{calc, bending, positioning}
\usetikzlibrary{arrows.meta}
\usetikzlibrary{shadings}

\usepgfplotslibrary{groupplots} 
\usetikzlibrary{pgfplots.groupplots} 
\usepgflibrary{shadings}
\usetikzlibrary{fadings,decorations.pathmorphing}
\pgfplotsset{compat=1.17}
\usetikzlibrary{shapes.geometric, arrows, positioning}
\usetikzlibrary{calc}
\usetikzlibrary{spy}
\usepackage{subcaption}

\usetikzlibrary{trees}
\usetikzlibrary{positioning,decorations.pathreplacing}
\usetikzlibrary{shapes.multipart}
\usetikzlibrary{chains,fit}
\usetikzlibrary{arrows.meta}
\usetikzlibrary{angles, quotes, babel}

\usetikzlibrary{patterns}
\usetikzlibrary{fadings}
\usetikzlibrary{positioning, arrows.meta, calc, matrix, fit, backgrounds}

\usepackage{longtable}
\usepackage{multirow} 

\tikzstyle{EDR}=[draw=cnrblue,preaction={clip, postaction={pattern=north west lines, pattern color=cnrblue}}]

\usepackage[colorlinks,allcolors={black},urlcolor={blue}]{hyperref} 

\usepackage{nicefrac}
\usepackage[ruled,vlined,linesnumbered,algo2e]{algorithm2e}

\begin{document}

\begin{frontmatter}

\title{Artificial neural network evaluation of geometric constants for polygonal domains\tnoteref{prep}
}

\author[bcpoli]{Beatrice Crippa}\ead{beatrice.crippa@polimi.it}
\affiliation[bcpoli]{organization={
             MOX-Laboratory for Modeling and Scientific Computing, Department of Mathematics, Politecnico di Milano},
             addressline={Via Bonardi 9},
             city={Milano},
             postcode={20133},
             country={Italy}}

\author[imati]{Sofia Imperatore\corref{cor1}}\ead{sofia.imperatore@imati.cnr.it}\cortext[cor1]{Corresponding author}
\author[imati]{Silvia Bertoluzza}\ead{silvia.bertoluzza@imati.cnr.it}
\author[imati]{Micol Pennacchio}\ead{micol.pennacchio@imati.cnr.it}

\affiliation[imati]{organization={Institute for Applied Mathematics and Information Technologies ``E. Magenes'', CNR },
                    addressline={Via Ferrata 5/A},
             city={Pavia},
             postcode={27100},
             country={Italy}}


\begin{abstract}
{
We propose an approach based on Artificial Neural Networks (ANNs) to evaluate geometric constants relevant to the analysis and design of numerical schemes for partial differential equations.
These constants play a central role, significantly influencing, for instance, a posteriori error estimates and the overall design of the computational strategy.
Our technique leverages ANNs to learn the dependencies between these constants and a set of descriptive geometric features associated to polytopal mesh elements. The main computational costs are confined to data processing and training phases, which can be performed offline once and for all.
This yields an {effective} tool for computing the constants, which we verify and show to be applicable across
different scenarios,
without substantial modifications -- demonstrating its broader usability beyond the specific example considered.}

\end{abstract}

\begin{keyword}

Geometric constants \sep deep learning \sep polygonal mesh \sep inequalities



\end{keyword}

\end{frontmatter}

\section{Introduction}\label{sec:intro}
A few key geometric constants—such as the Poincar\'{e} constant and those arising in inverse and trace inequalities—play a fundamental role in the analysis and design of numerical methods for solving partial differential equations (PDEs). Accurate estimation of these constants is crucial, as they directly influence both a priori and a posteriori error bounds, as well as the {design} of specific techniques, {such as stabilization methods including, for example, Nitsche’s method and SUPG.}
Reliable estimates enhance the robustness and efficiency of numerical schemes, particularly when dealing with complex geometries and high-order approximations.
{Despite their importance these constants are often employed either using overly pessimistic upper bounds or chosen in an almost arbitrary manner.}

While analytic expressions for these constants are generally unavailable, they can be obtained by solving a PDE eigenvalue problem. However, performing such a computation for each element in a tessellation would be prohibitively expensive. In the context of the Finite Element Method (FEM) -- where elements are typically triangles or {quadrilaterals} -- several studies have addressed the challenge of providing sharp estimates of these constants (see Section \ref{sec:known}). The situation becomes significantly more complex when dealing with polygonal tessellations.

Nevertheless, since the value of these constants depends on the geometry of the domain, one can attempt to relate them to a set of descriptive geometric features.
To this end, we propose leveraging Artificial Neural Networks (ANNs), a learning technique

recognized as a powerful tool for uncovering relationships between input data and numerical outputs. ANNs have proven highly effective when handling large datasets and are widely used for tasks such as classification, pattern recognition, prediction, and function approximation. In particular, ANNs excel at learning the underlying (possibly nonlinear) relationships between input and output data. In our case, we hypothesize that they can be trained to approximate these geometric constants associated with a polytope, given a set of descriptive geometric features.

This method is significantly less computationally expensive than directly computing the constants for each element of the mesh, as this computation is only required during the training dataset generation phase. This phase can be performed offline once and for all on a representative set of polytopes. Once the model is properly defined and trained, determining the constants for each mesh element is extremely fast.

In this work, we focus on three specific geometric constants: the Poincar\'{e} constant, and two additional constants arising from inverse and trace inequalities. We show
 that it is possible to obtain estimates for general polygonal domains, and that the structure of the proposed ANN-based approach 
 {remains largely unchanged when targeting different constants, indicating the generality and flexibility of the method beyond the specific cases considered.}

The paper is organized as follows.
In Section \ref{sec:known}, we specify the constants we aim to address and review their definitions along with some known results concerning their approximation. 
Sections \ref{sec:4} and \ref{sec:ANN} present methodological details that include selection of input attributes, 
tuning network hyperparameters, and strategies for optimizing the network architecture and improving the stability of the approach. Experimental results are reported in Section \ref{sec:5},
and finally, possible extensions and developments are outlined in Section \ref{sec:6}.

\section{Geometric constants}\label{sec:known}

The set of geometric constants relevant to the numerical analysis of PDEs
is quite extensive. It includes, among others, those appearing in the Poincar\'{e} inequality, trace and inverse inequalities, Sobolev embedding results, and  Korn’s inequality, to name a few. 
Most classical results concerning these constants have been developed within the
FEM framework, where—when focusing on the two-dimensional case—computational domains are typically triangles or rectangles. For such standard shapes, sharp estimates can often be derived analytically or computed efficiently.
An example of how some of these constants arise and can be computed in the context of FEM design and analysis can be found, for instance, in \cite{hughes_inequalities}.

Although these inequalities are well understood for simple geometries, the accurate estimation of their constants for general polygonal domains remains a significant challenge and a key motivation for the method proposed in this work.

Here, we focus on three key constants: the Poincar\'{e} constant, and those associated with trace and inverse inequalities,
showing the efficiency of the proposed ANN-based technique and that its structure remains largely unchanged when switching from one constant to another.
Throughout the rest of the paper, we assume $\el$ to be a non-degenerate polygon in $\mathbb{R}^2$.

\subsection{Poincar\'{e}–Friedrichs Inequality}

The Poincar\'{e}–Friedrichs inequality is a fundamental result in the theory of PDEs
providing a bound for a function in terms of its gradient. It plays a central role in analysis and numerical methods
by ensuring control of the \( L^2 \)-norm of a function through its first derivatives.

We then consider a classical Poincar\'{e}-type inequality of the form:
\begin{equation}\label{Poincare}
	\norm{w}_{0,\om}\leq \Cp \norm{\nabla w}_{0, \om} \qquad \forall \ w \in V,
\end{equation}
where 
\[
\norm{w}_{0, \om}^2 = \int_ \om | w |^2 \, dx, \qquad 
V \coloneqq \left\{w \in H^1(\om) \, : \, \int_{ \om} w \, dx = 0\right\}.
\]
It is well known that inequality \eqref{Poincare} holds with a constant $C_P$ that depends primarily on the geometry of $\el$.

Several analytical results exist for estimating \( \Cp \) in specific geometries. Exact values are known for rectangles, right triangles, and cubes \cite{exact_cp,lame,equilateral_cp}, while sharp two-sided bounds have been developed for triangles and tetrahedra \cite{eigenvalues}. \cite{isosceles_cp} derived upper bounds for isosceles triangles in terms of their diameter, and for a right isosceles triangle with unit-length sides, the exact constant is \( \Cp = 1/\pi \). Further bounds for convex polygons were presented in \cite{convex_cp}, showing that $
\Cp \leq \frac{d(\om)}{\pi},
$
where $d(\om)$ is the largest distance between any two vertices of the polygon.

\subsection{Trace Inequalities}\label{sec:trace}
The second inequality considered is a trace inequality, which provides bounds for a function on the boundary of a domain in terms of its volume norms. Such inequalities are critical in the analysis of partial differential equations and variational formulations, particularly when boundary conditions are imposed in a weak form \cite{brenner_scott,ern_guermond}. They also play a key role in domain decomposition methods \cite{toselli2005domain}, the analysis of boundary integral equations \cite{steinbach2008numerical}, stability estimates for non-conforming methods \cite{ainsworth2005posteriori}, and the derivation of both a priori and a posteriori error bounds \cite{verfurth2013posteriori}.
Moreover, trace inequalities are essential in the development of polytopal element methods (PEM), such as the Virtual Element Method \cite{beirao2013basic,gain2014virtual}, polygonal finite elements \cite{sukumar2004finite}, and mimetic finite difference methods \cite{da2014mimetic}, to name a few.

Here we focus on the following trace inequality, which holds for all \( u \in V \), where \( V = H^1(\om) \), and involves a constant \( \Ct \) that depends on the geometry of the domain:
\begin{equation}\label{trace_bound}
	\norm{u}^2_{0,\partial\om} \leq \Ct^2 \norm{u}_{0,\om} \left( h^{-1} \norm{u}_{0,\om} + \norm{\grad{u}}_{0,\om} \right) 
\end{equation}
where \( h \) denotes a characteristic length scale (e.g., the diameter of \( \om \)).

\subsection{Inverse Inequalities}\label{sec:inverse}
We finally consider inverse inequalities, where the expression inverse refers to the fact that, unlike the classical Poincaré inequality, they provide bounds on higher-order derivatives in terms of lower-order derivatives.

In particular, in this work we focus on the following inverse inequality, which holds for all polynomials \( p \in V=\mathbb{P}_1 \), i.e., polynomials of degree at most one defined on a polygon $\el$:
\begin{equation}\label{inverse_bound1}
	\norm{ \nabla p }_{0,\el}^2 \leq \Ci \, h^{-2}  \norm{ p }_{0,\el}^2,
\end{equation}
where 
\( h \) denotes the diameter of the polygon and 
$\Ci$ is a constant that, once again, depends on the shape of the domain 
$\el$, but not on its size.

These inequalities are 
particularly relevant in numerical analysis, especially in the context of
FEM and PEM,
where they characterize how discrete norms scale with respect to the local mesh size and are essential for stability analysis and error estimation.
Again, precise estimates are generally limited to standard reference elements (triangles, rectangles), and generalizations to polygonal domains are often empirical or require costly numerical evaluation \cite{ainsworth_inverse, repin, hughes_inequalities}.

\subsection{Computation of the Constants}\label{eigenvalue_problem}
We now recall how these constants can be computed 
showing that their evaluation can be reformulated as a generalized eigenvalue problem.

Let \( X \) and \( Y \) be Hilbert spaces defined on \( \el \) (or on its boundary), with \( V \subset X \cap Y \) a suitable subspace. For the inequalities described above, 
the goal is to estimate the best possible constant \( C_K \) in the abstract inequality:
\[
\norm{v}_{X(\el)} \leq C_K \norm{v}_{Y(\el)} \qquad \forall v \in V,
\]
where
\[
C_K= \sup_{v \in V} \frac{\norm{v}_{X(\el)}}{\norm{v}_{Y(\el)}}.
\]
{Let now $A$ and $B$ be the
operators induced by the inner products of 
$X$ and  $Y$, respectively. Then the constant $C_K$ can be equivalently expressed via the Rayleigh quotient:}
$$C_K^2 = \sup_{v \in V} \frac{(A v, v)}{(B v, v)}.$$
This leads to a generalized eigenvalue problem: find $\lambda \in \mathbb{R}$  and $u \in V $
such that
\begin{equation}
A u = \lambda B u.
	\label{eigen_gen}\end{equation}
The best constant is then given by:
$$C_K = \sqrt{\lambda_{\max}}, $$
where $\lambda_{\max}$ is the largest eigenvalue of the generalized eigenvalue problem associated with the operators $A$ and $B$.

However, computing such constants for each element in a mesh can be prohibitively expensive, especially in the case of complex or irregular polygonal meshes. While sharp estimates exist for standard elements such as triangles or quadrilaterals, general polygonal domains present significant challenges.
It is well established in the literature that the constants \( \Cp \), \( \Ct \), and \( \Ci \) are influenced primarily by geometric properties of the domain \( \el \). This motivates the idea of linking these constants to geometric descriptors of each polygon.

{To this end, we propose a deep learning approach based on ANNs to learn the
relationship between polygonal geometry and inequality constants.}

\section{Method description}\label{sec:4}
{In this Section, we provide details on the proposed method for the estimation of the geometric constants introduced above. The approach is based 
 on a \emph{supervised} learning regression scheme, where the model learns a continuous mapping between geometric descriptors of polytopal domains and the corresponding constants, by minimizing the error between predicted output and true constant values, computed via the eigenvalue formulation described in Section~\ref{eigenvalue_problem}.}

\begin{figure}[!t]
\centering
\begin{tikzpicture}[node distance=1.2cm and 0.8cm]

  \node (concave) {\includegraphics[width=0.18\linewidth]{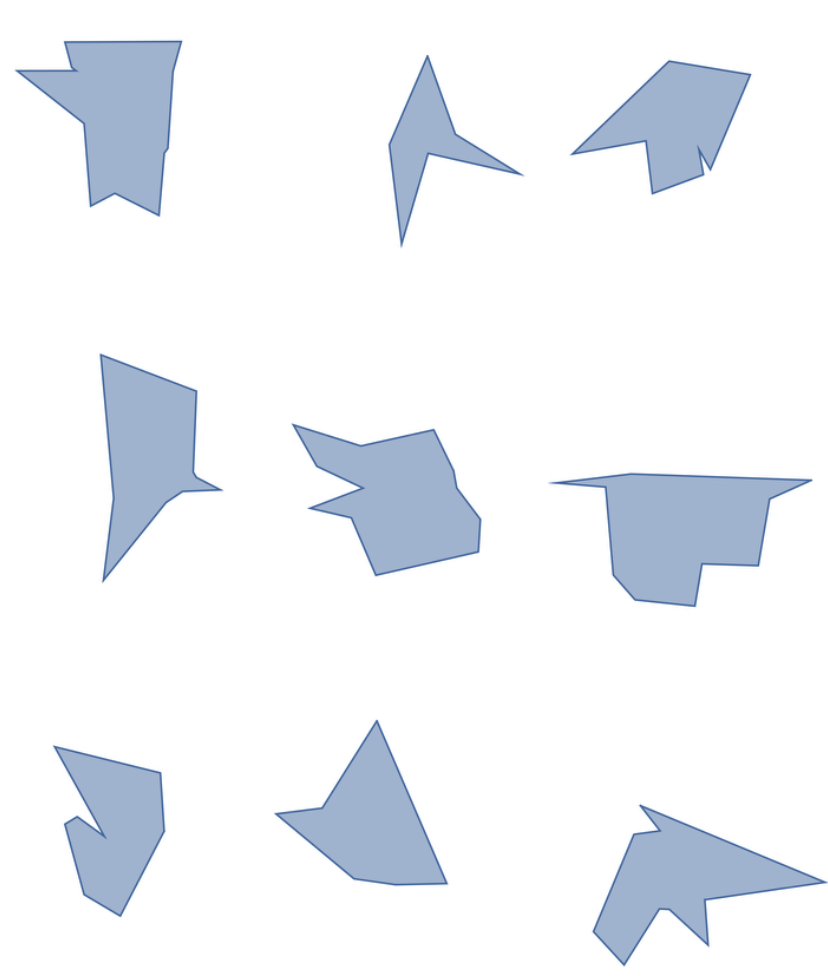}};
  \node (a) [below=2pt of concave] {\footnotesize (a) Non-convex polygons};

  \node[right= of concave] (mesh) {\includegraphics[width=0.18\linewidth]{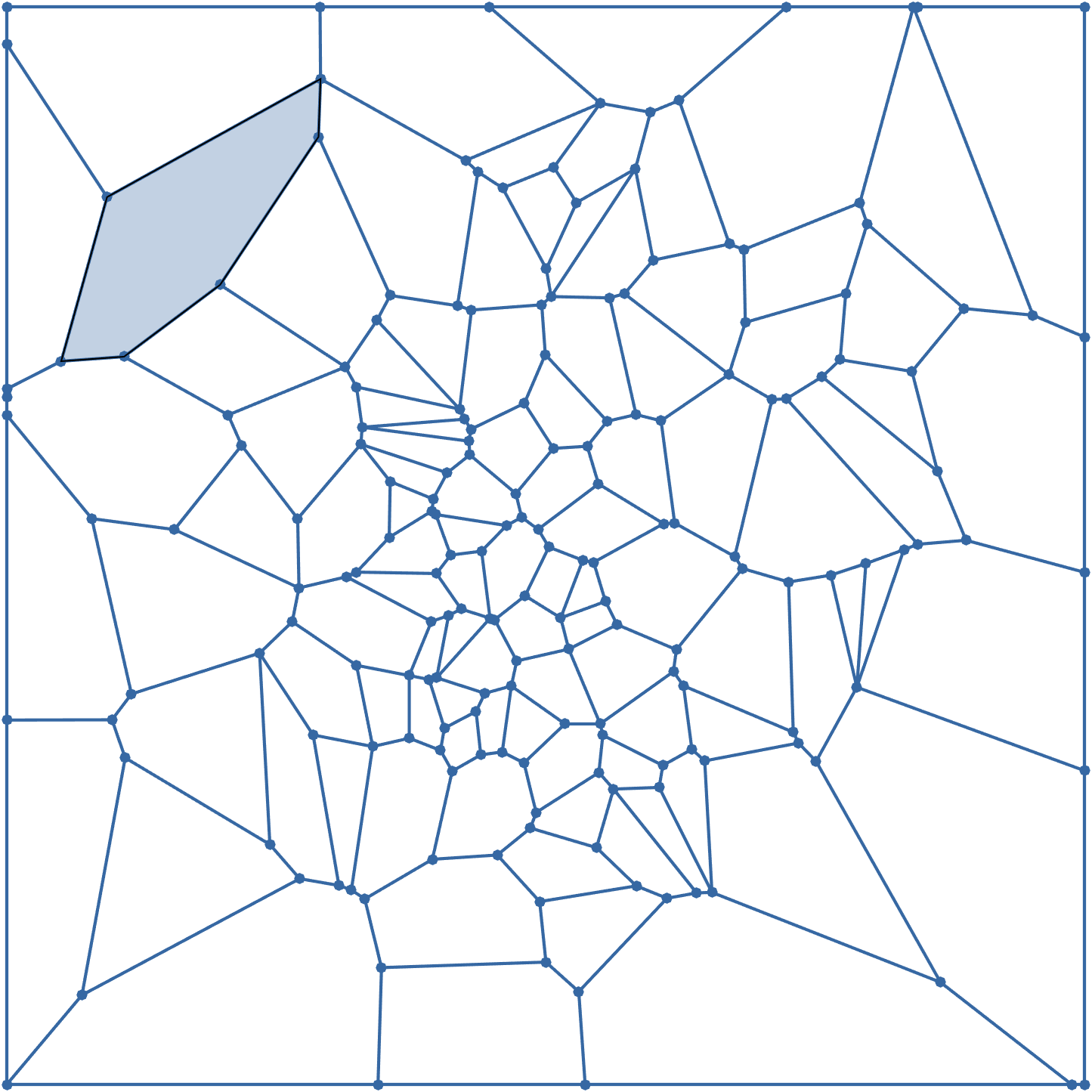}};

  \node[right=of mesh] (poly) {
    \begin{tikzpicture}[scale=10]
      \draw[fill=cnrblue!10, draw=cnrblue, thick]
        (0.0505, 0.6711) -- (0.1088, 0.6753) --
        (0.1980, 0.7420) -- (0.2893, 0.8794) --
        (0.2912, 0.9330) -- (0.0930, 0.8234) -- cycle;
    \end{tikzpicture}
  };
  \draw[-{Latex[round]}, thick] 
    ($(mesh.north east)+(-2.1cm,-.7cm)$) -- ($(poly.west)+(0cm,0.25cm)$);

  \node (b) at ($(mesh)!0.5!(poly)$) {};

    \path
    	let 
    	\p1 = (a),
    	\p2 = (b),
    in  node at (\x2,\y1) {\footnotesize(b) Voronoi mesh and extracted convex polygon};
 \end{tikzpicture}
\caption{Example of non-convex polygons (a); a Voronoi mesh and convex polygon extracted from the mesh (b).}\label{fig:polygons}
\end{figure}

More precisely, our study focuses on general planar domain (either non-convex or convex), characterized by polygons with an arbitrary number of edges.
Figure~\ref{fig:polygons}(a) shows an example featuring nine randomly generated non-convex polygons obtained  using the procedure described in Appendix A (Algorithm~\ref{alg:nonconvex-gen}). In addition to these, we also investigate the subspace of convex polygons to evaluate the potential benefits that this geometric constraint may offer to the overall procedure.
In the convex case, the polygons are extracted from a Voronoi mesh generated over a square domain in $\mathbb{R}^2$, such as 
$[0,1]^2$.
An illustrative example is provided in Figure~\ref{fig:polygons}(b).

We also recall that the constants introduced in Section~\ref{sec:known} scale with the diameter of the  domain. For example, the Poincaré constant satisfies:
\[
C^P_\om = \text{diam}(\om) \cdot \widetilde{C}^P_\om,
\]
where $\widetilde{C}^P_\om$ is the constant associated with the rescaled domain $\widetilde{\om}$, obtained by normalizing $\om$ to have unit diameter.
Therefore, we restrict our analysis to polygons with unit diameter. The generalization to polygonal domains with arbitrary diameter is discussed in Section~\ref{general_polygons}.

Since the constants we want to evaluate, 
$C_p$, $C_i$, and $C_t$, are known to depend on the polygon $\domain$ 
we aim to construct a continuous map 
\begin{equation}\label{eq:map}
    f:\domain \longmapsto C_\domain,
\end{equation}
which associates each domain with its corresponding constant.

Polygonal domains $\domain$ 
can be described in multiple ways, such as by specifying the coordinates of their vertices or through image-based representations. The constants introduced earlier are domain-dependent, and we 
postulate that they can be expressed in terms of a set of geometric features. More specifically, following 
the approach of \cite{attene2021}, we represent each polygon using a vector of \emph{geometric quality quantities}, denoted as 
$$\mathbf{m}(K) = \{ m_1(K),\dots, m_d(K) )\}, $$ which we refer to simply as metrics. Each component $m_i(K)$  corresponds to a geometric feature, as detailed in Table~\ref{tab:quality_metrics}.

\begin{table}[!t]
  \centering
  \caption{Geometric quality metrics describing each polygon \(K\).}
  \label{tab:quality_metrics}

\resizebox{.9\textwidth}{!}{
\begin{minipage}{0.5\linewidth}
   \begin{tabular}{@{}l p{3.cm} p{4.5cm}  @{}}
     \toprule
        \textbf{Name} && \textbf{Description} \\
     \midrule
     \textbf{CC}     & Circumscribed circle radius & Radius of the smallest circle enclosing \(K\), computed via Welzl’s algorithm.  \\
     \textbf{IC}     & Inscribed circle radius     & Radius of the largest circle fully contained in \(K\). \\
     \textbf{CR}     & Circle ratio                & \(\mathrm{IC}/\mathrm{CC}\in[0,1]\); indicates circularity. \\
     \textbf{AR}     & Area                        & Total area of polygon \(K\). \\
     \textbf{KE}     & Kernel area                 & Area of the kernel (visibility region); equals AR if \(K\) is convex. \\
     \textbf{KAR}    & Kernel-area ratio           & \(\mathrm{KE}/\mathrm{AR}\in[0,1]\); fraction of area that is kernel. \\
     \textbf{APR}    & Area–perimeter² ratio       & \(\mathrm{AR}/(2\rho)^2\); area scaled by squared perimeter. \\
     \textbf{SE}     & Shortest edge               & Length of the shortest edge in \(K\). \\
     \textbf{sSE}    & Scaled shortest edge        & \(\mathrm{SE}/\mathrm{CC}\in[0,1]\); normalized by CC. \\
     \textbf{ER}     & Edge ratio                  & \(\mathrm{SE}/\mathrm{LE}\in(0,1]\); ratio of shortest to longest edge.\\
     \textbf{MPD}    & Min. point-to-point dist.   & Minimum distance between any two vertices of \(K\). \\
    \textbf{sMPD}   & Scaled min. point to point  dist.       & \(\mathrm{MPD}/\mathrm{CC}\in[0,1]\); normalized by CC. \\
     \textbf{MA}     & Minimum inner angle         & Smallest internal angle in \(K\).\\
     \textbf{MX}     & Maximum inner angle         & Largest internal angle in \(K\). \\
     \bottomrule
   \end{tabular}
\end{minipage}\begin{minipage}{0.5\linewidth}
  \hspace*{3cm}
  \centering
  \begin{tikzpicture}[scale=0.7]  
		
 		\def\rKE{1.0}
 		\def\rIC{1.8}
 		\def\rCC{3.5}
 		\def\nSpikes{10}
 		\def\angleStep{360/\nSpikes}
		
 		\coordinate (O) at (0,0);
		
 			\coordinate (O) at (0,0);
		
 		\fill[red!60] (O) circle (\rKE);
 		\node at (O) {\textbf{KE}};

 		\foreach \i in {0,...,9} {
 			\pgfmathsetmacro{\aStart}{\i * \angleStep}
 			\pgfmathsetmacro{\aEnd}{(\i + 1) * \angleStep}
 			\pgfmathsetmacro{\aMid}{(\aStart + \aEnd)/2}
			
 			\path[fill=gray!60]
 			(\aStart:\rIC) -- (\aMid:\rCC) -- (\aEnd:\rIC) -- cycle;
 		}

			\fill[gray!20] (O) circle (\rIC);
			
 			\fill[cnrblue!60] (O) circle (\rKE);
 			\node at (O) {\textbf{KE}};
			
 			\draw[cnrblue,dashed,thick] (O) circle (\rIC);
 			\node at (90:\rIC+0.3) {\textbf{IC}};
 			\draw[cnrblue,dashed,thick] (O) circle (\rCC);
 			\node at (110:\rCC+0.4) {\textbf{CC}};

 	\path (36:\rIC) coordinate (P1);
 	\path (72:\rIC) coordinate (P2);
 	\path (54:\rCC) coordinate (P2CC);
 	\draw[decorate,decoration={brace,amplitude=6pt},thick]
 	(P1) -- (P2) node[midway,right=5pt] {\textbf{MPD}};
	\fill (P1) circle (2pt);
 	\fill (P2) circle (2pt);
 	\node at ($(P1)!0.6!(P2)+(0.3,1.3)$) {\textbf{SE}};
 		\fill (P2CC) circle (2pt);
  \draw[line width=1pt] (P2) -- (P2CC);
	
 		\node at (0,-1.3) {\textbf{AR}};
				
 	\end{tikzpicture}

 \end{minipage}
 }
   \end{table}

This representation is particularly advantageous in three-dimensional settings, as it enables the characterization of a polytope without relying on its exact—often computationally intensive—definition. Instead, it leverages fixed, measurable quantities intrinsic to the polytope.
Therefore, we postulate that,  for each polygon $K$,  the feature vector $\mathbf{m}(K)$ contains sufficient information to determine the associated constant $C_K$.
To model this relationship, we define the mapping
\begin{equation}\label{eq:mapann}
    f: \text{\textbf{m}(\domain}) \longmapsto C_\domain
\end{equation} 
using a compact feed-forward neural network, specifically a multi-layer perceptron (MLP), which captures the underlying dependency between the computational domain and the target constant.
In practice, the metrics $\mathbf{m}$ serve as input to the MLP, which outputs an estimate of the corresponding constant $C_\domain(\mathbf{m})$. 
A visual overview of this methodology is provided in Figure~\ref{fig_ANN}.

\begin{figure}[!t]
		\centering

            \begin{tikzpicture}[
    node distance=1.5cm and 1.5cm,
    every node/.style={align=center},
    arrow/.style={black, -{Latex[round]}, thick},
    rescale/.style={black, font=\scriptsize},
    neuron/.style={circle, draw, minimum size=6pt, inner sep=0pt},
    layer/.style={matrix of nodes, nodes=neuron, column sep=4mm, row sep=4mm, anchor=center},
    stage/.style={draw=cnrblue!50, fill=cnrblue!10, thick, rounded corners, inner sep=0.5em}
  ]
  
  \node (X) {$K$};
  
  \node[right=of X] (Xt) {$\mathbf{m}$};

  \draw[arrow] (X) -- (Xt);
  
  \node (scale1) [rescale, below=0.1cm of $(X)!0.5!(Xt)$] { }; 

  \node (output-1-1) [right = of Xt] {
    \begin{tikzpicture}[scale=0.5]
					\tikzstyle{unit}=[draw, shape=circle, minimum size=0.035cm]
					
					\tikzstyle{hidden}=[draw,shape=circle, minimum size=0.035cm]
					
					\node[unit](x0) at (0,3.5){};
					\node at (0,2.8){$\vdots$};
					\node[unit](x1) at (0, 1.75){};
					\node at (0,1){$\vdots$};
					\node[unit](xd) at (0,0){};
					
					\node[hidden](h10) at (3-1.5,4){};
					\node[hidden](h11) at (3-1.5,2.5){};
					\node at (3-1.5,1.5){$\vdots$};
					\node[hidden](h1m) at (3-1.5,-0.5){};
					
					\node(h22) at (5-1.5,0){};
					\node(h21) at (5-1.5,2){};
					\node(h20) at (5-1.5,4){};
					
					\node(d3) at (6-2,0){$\dots$};
					\node(d2) at (6-2,2){$\dots$};
					\node(d1) at (6-2,4){$\dots$};
					
					\node(hL12) at (7-1.5,0){};
					\node(hL11) at (7-1.5,2){};
					\node(hL10) at (7-1.5,4){};
					
					\node[hidden](hL0) at (9-1.5,4){};
					\node[hidden](hL1) at (9-1.5,2.5){};
					\node at (9-1.5,1.5){$\vdots$};
					\node[hidden](hLm) at (9-1.5,-0.5){};

					\node[unit](y2) at (11-2,2){};

					\draw[-latex] (x0) -- (h11);
					\draw[-latex] (x0) -- (h1m);
					
					\draw[-latex] (x1) -- (h11);
					\draw[-latex] (x1) -- (h1m);
					
					\draw[-latex] (xd) -- (h11);
					\draw[-latex] (xd) -- (h1m);

					\draw[-latex] (hL0) -- (y2);

					\draw[-latex] (hL1) -- (y2);

					\draw[-latex] (hLm) -- (y2);

					\draw[-latex, path fading=east] (h10) -- (h21) node[pos=0.3, sloped, above]  {};
					\draw[-latex, path fading=east] (h10) -- (h22);
					
					\draw[-latex, path fading=east] (h11) -- (h21);
					\draw[-latex, path fading=east] (h11) -- (h22);
					
					\draw[-latex, path fading=east] (h1m) -- (h21);
					\draw[-latex, path fading=east] (h1m) -- (h22);
					
					\draw[-latex, path fading=west] (hL10) -- (hL1);
					\draw[-latex, path fading=west] (hL11) -- (hL1);
					\draw[-latex, path fading=west] (hL12) -- (hL1);
					
					\draw[-latex, path fading=west] (hL10) -- (hLm);
					\draw[-latex, path fading=west] (hL11) -- (hLm);
					\draw[-latex, path fading=west] (hL12) -- (hLm);

				\end{tikzpicture}
    
    };

  \node[right=1.2cm of output-1-1] (yt) {$C_K\left(\mathbf{m}\right)$}; 
%
  \draw[arrow] (Xt) -- (output-1-1);
  \draw[arrow] (output-1-1) -- (yt);
  \node (scale2) [rescale, below=0.1cm of $(yt)$] { }; 
  
\begin{scope}[on background layer]
  \node[stage, fit=(X)(Xt)(scale1)] (inputbox) {};
  \node[above=0.2cm of inputbox] {Input Features};

  \node[stage, fit=(output-1-1)] (nnbox) {};
  \node[above=0.2cm of nnbox] {Neural Network};

  \node[stage, fit=(yt)(scale2)] (outbox) {};     
  \node[above=0.2cm of outbox] {Estimate};  
\end{scope}

\end{tikzpicture}

\caption{Overview of the proposed approach: the geometric metrics are computed for a polygon $K$ 
        to obtain the evaluation of the considered constant corresponding to $K$.}\label{fig_ANN}
\end{figure}

We then filter the metrics to retain only those that are meaningful and non-redundant for each polygon class considered.
For convex polygons, certain metrics become redundant: for instance, the kernel area equals the total area, the kernel-area ratio is always 1, and the shape regularity coincides with the circle ratio.
These metrics are therefore excluded, and the selected sets for both polygon types are summarized in Table~\ref{metric_polygons}.

\begin{table}[!t]
  \centering
  \caption{Metrics for convex and non-convex polygons.}
  \begin{tabular}{@{}ll@{}}
    \toprule
    \textbf{Polygon Type} & \textbf{Metric \(\mathbf{m}\)} \\
    \midrule
        Non-convex polygon &
    \(\{\CR,\ \AR,\ \KE,\ \KAR,\ \APR,\ \sSE,\ \ER,\ \sMPD,\ \MA,\ \MX,\ \SR,\ \ISO\}\) \\
    Convex polygon &
    \(\{\CR,\ \AR,\ \APR,\ \sSE,\ \ER,\ \sMPD,\ \MA,\ \MX,\ \ISO\}\) \\
    \bottomrule
  \end{tabular}\label{metric_polygons}
\end{table}

\section{Data-Driven geometric constant evaluation via Artificial Neural Networks}\label{sec:ANN}

In this section we present the proposed deep learning approach based on ANNs to learn the underlying mapping from polygon geometry to inequality constants.  
By training on a dataset of polygons
and their associated constants (computed via the eigenvalue formulation), the model can provide accurate and computationally efficient predictions for any polygonal domains. This approach enables practical applications in PEMs, where such constants have to be evaluated efficiently
on a wide range of domain geometries.

\subsection{Dataset generation}\label{sec:datagen}
{For both non-convex and convex class, we generate polygons
characterized by unitary diameter; subsequently, for each polygon we compute the metrics defined into Table \ref{metric_polygons} and subdivede the data into training and validation set and testing set.}
Finally, we eliminate outliers from the generated data sets.
More precisely, we consider outliers the samples with at least one attribute that has z-score $\geq 2$.
Figure~\ref{fig:boxplot} shows the boxplots of the attributes in the training set, before (left) and after (right) the outlier elimination both for non-convex (top) and convex (bottom) training sets. This pre-processing gives a dataset with a better distribution around its mean value, hence it has a regularizing effect on the data-driven method.
Finally, by letting $M$ be the total amount of polygons, the generated dataset consist of label data
    \begin{equation}\label{eq:dataset}
        \left\{\left(\mathbf{m}_i, C_{K_i}\right)\ :\  i = 1, \ldots, M\right\}.
    \end{equation}

    \begin{figure}[!t]
    \centering
    \include{boxplot_nonconvex}\vspace*{-2cm}
    \include{boxplot_convex}
    \caption{Boxplots of the attributes in the non-convex (top) and convex (bottom) training sets, before (left) and after (right) outliers elimination.}
    \label{fig:boxplot}
\end{figure}

\subsection{ANN architecture}\label{sec:architecture}
The assumption of continuity of the map \eqref{eq:mapann} implies that it can be approximated by a feedforward ANN taking a description of the domain as input and giving the value of the corresponding constant as output.
ANNs are computational learning systems made of interconnected nodes, called neurons or units, grouped into layers, that transom input data into numeric output, learning the feature relation during the so-called training phase.
Different {ANN} architectures determine different learning models, where the word {architecture} refers to the overall structure of the network: how many units it has and how these units are connected to each other.

To approximate $f$, we consider a MLP
consisting of one input layer, $L\in \mathbb{N}_{\geq1}$ fully connected hidden layers, and one output layer; all pair of consecutive layers are alternated by non-linear activation functions. 
More precisely,
let $L\in \mathbb{N}_{\geq1}$, where $L+2$ is the total number of layers of a {MLP} architecture. For each $\ell=0, \ldots, L+1$, let $N_\ell$ be the number of hidden units for the $\ell$-th layer and let
$
\sigma_\ell:\mathbb{R}\to\mathbb{R}
$ be a non-linear activation function. 
In the input layer, each unit corresponds to a feature of the
geometric quality metrics $\mathbf{m}$, described in Table~\ref{tab:quality_metrics}.
Hence, there are as many units $N_0$ as the number of features of the considered data, thus in our specific case $N_0 = 12$ for non-convex polygons, whereas $N_0=9$ for convex ones, as illustrated in Table~\ref{metric_polygons}. Moreover, the direct graph connections are only in output, from the units of the input layers towards the units of the following one. Subsequently, $L$ hidden layers with $N_\ell$ units for $\ell=1, \ldots, L$, are placed between the input and output layers: at any hidden layer, each unit receives the same real value from every unit of the previous layer and produces a real value that is transomed via the activation function $\sigma$ and passed to every unit of the following layer.
Finally, the output layer has many units $N_{L+1}$ as the function at hand, which the {ANN} is required to approximate, hence in our case $N_{L+1} = 1$. There is no other following layer, and the output unit connections are only along the incoming direction.
Let $\omega_{ij}^{(\ell)}$ be the weight associated with the edge between the $i$-th unit of the $\ell$-th layer and the $j$-th unit of the $(\ell-1)$-th layer, and let $b_{i}^{(\ell)}$ be the bias associated with the $i$-th unit of the $\ell$-th layer, for $\ell=0, \ldots, L+1$. Then, the output of each neuron of the MLP can be expressed as
\[
    a_{j}^{(\ell)} = \sigma_\ell\left(\sum_{k=1}^{N_{\ell-1}}\omega_{ij}^{(\ell)}a_k^{(\ell-1)} + b_j^{(\ell)}\right),
\]
for $j=1, \ldots, N_\ell$ and $\ell = 0, \ldots, L+1$.
The optimal values of the weights $\omega_{ij}^{(\ell)}$ for $i,j=1, \ldots, N_\ell$ and $\ell = 0, \ldots, L+1$ are set during the training of the MLP model.

The choice of activation functions is a critical part of neural network design, and it depends on the problem at hand. In general,
they are needed to be non-linear functions in order to avoid the trivial situation in which the output is only a linear combination of the input data \cite{goodfellow2016}. In addition, their concurrence allows to represent arbitrarily complex functions \cite{aggarwal2018}.
More precisely for its simple piecewise linear structure and high performance \cite{yaro2017}, we choose $\sigma_\ell = \sigma$ for $\ell=0, \ldots, N+1$ to be the {ReLU} activation function, namely
\[
    \sigma(x) = \max\{0, x\},\ \text{for each}\ x\in\mathbb{R}.
\]

The overall number of layers $L+2$ gives the \emph{depth} of the model.
The total number of learnable parameters $\#LP$ of a {MLP} architecture can be computed as,
     \[
	\#LP \coloneqq \sum_{\ell=1}^{L+1}N_\ell\left(N_{\ell-1}+1\right).
    \]
By adding more layers and more units within a layer, a deep {ANN} can represent functions of increasing complexity, but of course the computational costs increase, as well as the risk of overfitting. It is fundamental to find a proper architecture that is a trade-off between computational training costs and the accuracy that can be achieved.

\subsection{Loss function}
 As introduced at the beginning of Section~\ref{sec:4}, the proposed method is based on a supervised learning regression scheme performed by a MLP.
 To optimize the parameters characterizing the MLP, first of all 
 we divide the dataset, generated and processed as described in Section~\ref{sec:datagen}, into a training set T and a validation set V, taken as 70\% and 30\% of the data, respectively, such that $T\cap V=\emptyset$. The model will be trained on T and at the same time validated on V.
 More precisely, the training process consists in minimizing the loss function $\mathcal{L}$ defined as the
 mean squared error between the output of the regression model and the true target value.
 In particular, 
 \begin{equation}\label{eq:loss}
     \mathcal{L}\left(\left\{\mathbf{m}_i, C_{K_i} \right\}_{i=1}^M; W, \mathbf{b}\right) = \frac{1}{2m}\sum_{i=1}^M\left\|C_{K_i}\left(\mathbf{m_i}\right) - C_{K_i}\right\|_2^2,
 \end{equation}
 where $W$ is the collection of the model weights $\omega^{(\ell)}_{ij}$, $\mathbf{b}$ the collection of biases $b^{(\ell)}_i$, for $i,j=1, \ldots, N_{\ell}$, and $\ell = 0, \ldots, L+1$, and $\left\{\mathbf{m}_i, C_{K_i} \right\}_{i=1}^M$ are the training data contained in T.
 At the end of each training epoch, the model is validated by evaluating \eqref{eq:loss} on the validation dataset V.

\subsection{Training and hyperparameter selection}
As described in Section~\ref{sec:architecture}, we consider a dense MLP, characterized by $L$ hidden layers, consisting of $N_\ell = N$ units, for each $\ell = 1, \ldots, L$, and with ReLU activation functions after each layer.
To solve the training optimization problem \eqref{eq:loss}, we employ the Adam optimizer, with learning rate $\eta$.
In order to optimize the performance of the model, we set up an extensive grid search 
on the depth $L$, the width $N$ and the learning rate $\eta$, thus we compare the model performance across $1500$ training epochs on the validation set V for different combinations of these three values. More precisely, for each considered constant, for non-convex and convex polygons, we evaluate the accuracy of the model in terms of validation loss for all the hyperparameter values, summarized in Table~\ref{tab:hyper}.
We consider as best models the one with minimum validation loss.
Figure~\ref{fig:manual-tuning} illustrates all validation losses along 1500 epochs, for different values of $L$ from left to right, $N$, and $\eta$ from top to bottom for the hyperparameter tuning of the MLP predicting the Poincar\'{e} constant on non-convex polygons.
\begin{table}[!t]
    \centering

\begin{tabular}{ll}
\toprule
\# 
experiments & 160 \\
Tested $L$ parameters & $1, 3, 5, 7$ \\
Tested $N$ parameters & $8, 16, 32, 64, 128, 256, 512, 1024$ \\
Tested learning rates & $10^{-4}, 5\cdot10^{-4}, 10^{-3}, 5\cdot10^{-3}, 10^{-2}$ \\
\bottomrule
\end{tabular}
    \caption{Summary of the hyperparameter values considered to finetune the MLP architecture: the depth $L$, the width $N$ and the optimizer learning rate $\eta$.}
    \label{tab:hyper}
\end{table}

\begin{figure}[!t]
    \centering
    \input{nonconvex_poincare_tikz_tuning}
    
 \caption{Non-convex polygons -- Poincar\'{e} ineqiality. Validation loss evaluated at each of the first $1500$ training epochs for various values of network depth $L \in \{1, 3, 5, 7\}$, width $N \in \{8, 32, 64, 128, 256, 512, 1024\}$, and learning rate $\eta \in \{10^{-2},\ 5{\cdot}10^{-3},\ 10^{-3},\ 5{\cdot}10^{-4},\ 10^{-4}\}$.}
    \label{fig:manual-tuning}
\end{figure}

\section{Results}\label{sec:5}
In this section, we present the ANN-based estimation of geometric constants for the Poincar\'{e}, inverse and trace inequality introduced in Section~\ref{sec:known}, on non-convex and convex domains.
{For the non-convex and convex scenarios, and for each inequality considered, we generate training, validation and test sets as described in Section~\ref{sec:datagen}. More precisely, we generate dataset as \eqref{eq:dataset}, with $M = 500$ for training and validation, and $M = 75$ for testing. Each dataset is subsequently further processed to remove outliers, as explained in Section~\ref{sec:datagen}, 
resulting in $383$ samples for training and validation, and $59$ samples for testing in for non-convex polygons; on the other hand, for convex polygons, we have $385$ samples for training and validation, and $53$ samples for testing.}

It is well known that the validation loss is consistent with the scale of the target values. Therefore, we expect and have indeed verified differences in the best validation results across datasets, as the magnitude of the computed constants varies significantly. For this reason, for each training dataset considered, we display the mean, standard deviation, minimum, and maximum values of the target constants contained in each training set~\eqref{eq:dataset}.
These standard statistical measures are reported to characterize the data distributions and support comparisons across different datasets.
{In the final part of this section, we illustrate how to extend the proposed method to polygons of arbitrary diameter.}

The implementation of the learning method is based on the TensorFlow Python library \cite{tensorflow}; geometric features of each polygon are computed using the scikit-geometry package \cite{WinNT}
while the approximation of the constants in the training set—obtained via the generalized eigenvalue problems described in Section \ref{eigenvalue_problem}—is performed using the linear FEM, implemented in MATLAB R2024a, on a mesh with element diameter equal to $1/40$ of the polygon diameter.

\subsection{Poincar\'{e} inequality on non-convex polygons}
We start by considering the Poincar\'{e} inequality and estimate the corresponding geometric constant on non-convex polygons.
Figure~\ref{fig:nonconvex_poincare}~(a) displays the mean, standard deviation, minimum, and maximum values for the target Poincar\'{e} constant of the training data set.
The training data shows a mean of 0.2458 (std = 0.0270), ranging from 0.2027 to 0.3323. 
Such a small standard deviation, relative to the mean, indicates a low degree of variability in the dataset, thus this limited spread suggests that most observations lie close to the average, with no substantial deviations.

The selection of hyperparameters leads to an MLP architecture with 1 hidden layer and 128 units to be trained with the Adam optimizer with learning rate $10^{-2}$.
The training and validation losses of this model over $1500$ epochs are illustrated in Figure~\ref{fig:nonconvex_poincare}~(b), in particular
the minimum validation loss with value $3.5457\cdot10^{-4}$ is achieved at epoch 514; subsequently the model tends to overfit the training data and the validation loss to stagnates. 
Finally, this trained model achieves a test error of $2.3776\cdot10^{-4}$ in the test data set; therefore, the proposed model is able to predict the geometric constant related to the Poincar\'{e} inequality on non-convex polygons with high accuracy.
The hyperparameters and error values for this experiment are summarized in Table~\ref{tab:results_all}.

\begin{figure}[!t]
    \begin{subfigure}[b]{\linewidth}
            \centering
\begin{tikzpicture}
\begin{axis}[enlargelimits=true,
width=\linewidth,
height=.45\linewidth,
ylabel={constant value},
xlabel={\#Non-convex polygon data - Poincar\'{e} inequality},
grid=major, grid style={dashed,gray!30},
grid style=dashed,
xmin = 0,
xmax = 385,
ymin=.2,
ymax=0.4,
legend style={at={(rel axis cs:0.5,-0.25)}, legend columns=6,fill=none,draw=none,anchor=north, /tikz/every even column/.append style={column sep=12pt}},
]

\addplot+[cnrblue, every mark/.append style = {fill = cnrblue}, mark = *, mark size = 1.5pt, mark options = {solid}, only marks] table [x=, y=y, col sep=comma]{hypertuning/clean_train_nonconvex_poincare_labels.csv};\addlegendentry{data}

\addplot+[color=red, thick, mark size = 0pt] coordinates {(-50,0.247293) (500,0.247293)};\addlegendentry{mean}

\addplot+[color=dgreen, thick, densely dashed, mark size = 0pt] coordinates {(-50,0.278428) (500,0.278428)};\addlegendentry{$+$std}

\addplot+[color=olive, thick, densely dashed, mark size = 0pt] coordinates {(-50,0.2161580) (500,0.2161580)};\addlegendentry{$-$std}

\addplot+[color=cyan, thick, densely dashed, mark size = 0pt] coordinates {(-50,0.195922) (500,0.195922)};\addlegendentry{min}

\addplot+[color=magenta, thick, densely dashed, mark size = 0pt] coordinates {(-50,0.427265) (500,0.427265)};\addlegendentry{max}

\end{axis}
\end{tikzpicture}
            \caption{}
    \end{subfigure}

    \begin{subfigure}[b]{\linewidth}
            \centering
\begin{tikzpicture}
\begin{axis}[enlargelimits=true, ymode = log,
width=\linewidth,
height=.3\linewidth,
xlabel={epochs},
ylabel={MSE},
grid=major, grid style={dashed,gray!30},
grid style=dashed,
xmin=0, xmax=1500,
ymax=0.1,
ymin=0.00001,
ytick={1e-1, 1e-2, 1e-3, 1e-4, 1e-5},
yticklabels={$10^{-1}$, $10^{-2}$, $10^{-3}$, $10^{-4}$, $10^{-5}$}]

\addplot+[cnrblue, every mark/.append style = {fill = cnrblue}, mark = *, mark size = 0.5pt, mark options = {solid}, opacity = 0.65] table [x=, y=trainloss, col sep=comma]{hypertuning/nonconvex/poincare/hypertuning_nonconvex_poincare_L_1_N_128_lr_0.01_losses.csv};\addlegendentry{training}

\addplot+[orange, every mark/.append style = {fill = orange}, mark size = 0.5pt, mark options = {solid}, opacity = 0.35] table [x=, y=valloss, col sep=comma]{hypertuning/nonconvex/poincare/hypertuning_nonconvex_poincare_L_1_N_128_lr_0.01_losses.csv};\addlegendentry{validation}
\end{axis}
\end{tikzpicture}
    \caption{}
    \end{subfigure}
\caption{(a) Statistical analysis of the training target labels for learning the Poincar\'{e} constant on non-convex polygons: the data (blue), the mean value (read), the standard deviation (green and olive), the minimum (cyan) and the maximum values (magenta). (b) Training and validation loss over epochs for Poincar\'{e} inequality on non-convex polygons.}\label{fig:nonconvex_poincare}
\end{figure}

\begin{table}[!bh]
    \centering
    \footnotesize
    \setlength{\tabcolsep}{4pt}
    \resizebox{\linewidth}{!}{
    \begin{tabular}{r c cc c cc c cc}
        \toprule
         && \multicolumn{2}{c}{Poincar\'{e}} && \multicolumn{2}{c}{Inverse} && \multicolumn{2}{c}{Trace} \\
         && non-convex & convex && non-convex & convex && non-convex & convex\\
        \cmidrule{3-3}\cmidrule{4-4} \cmidrule{6-6}\cmidrule{7-7} \cmidrule{9-9}\cmidrule{10-10}
        \# layers && $1$ & $1$ && $7$ & $3$ && $3$ & $7$\\
        \# units && $128$ & $64$ && $1024$ & $256$ && $32$ & $32$\\
        $\eta$ && $10^{-2}$ & $5\cdot10^{-3}$ &&  $10^{-2}$ & $10^{-3}$ && $5\cdot10^{-3}$ & $5\cdot10^{-3}$\\ 
        Epoch && $514$ & $1284$ && $1139$ & $548$ && $1209$ & $1352$\\
        Val. loss && $3.5457  \cdot 10^{-4}$ & $2.3490 \cdot 10^{-5}$ && $4.1046  \cdot 10^{-1}$ & $3.0341 \cdot 10^{-2}$ && $5.7048  \cdot 10^{-2}$ & $3.1308  \cdot 10^{-2}$\\
        Test error && $2.3776 \cdot 10^{-4}$ & $3.7580 \cdot 10^{-5}$ && $5.1935 \cdot 10^{-1}$ & $9.1029 \cdot 10^{-1}$ && $5.5408 \cdot 10^{-2}$ & $7.5794 \cdot 10^{-2}$\\
        \bottomrule
    \end{tabular}
    }
    \caption{Summary of results for Poincar\'e, Inverse and Trace inequality both on non-convex and convex polygons. For each configuration we report the hyperparameter of the MLP architecture: number of hidden layers (\# layers), the number of units per hidden layer (\# units), and the learning rate of the Adam optimizer ($\eta$); the epoch ast which the model achieves the best validation loss, the best validation loss, and the corresponding test error registered by the selected model.}\label{tab:results_all}
\end{table}

\subsection{Poincar\'e inequality on convex polygons}
For convex polygons, the training data shows a slightly narrower spread, \ie the training data exhibit even lower variability and are more concentrated near the average. More precisely, the data are characterized by a mean of $0.2217$
(std = $0.0136$), ranging from $0.2217$ to $0.2968$.
Figure~\ref{fig:convex_poincare}~(a) displays the mean, standard deviation, minimum, and maximum values for the target Poincar\'{e} constant of the training data set.

In this case, the hyperparameter selection leads to an MLP architecture with 1 hidden layer, as in the previous experiment, but with less unit per layer, namely 64; moreover, the learning rate of the Adam optimizer is $5\cdot10^{-3}$.
The training and validation losses of this model over $1500$ epochs are illustrated in Figure~\ref{fig:convex_poincare}~(b);
{where we can observe that, differently from the non-convex configuration, the model does not tend to overfit. More preciesly,}
the minimum validation loss \eqref{eq:loss} with value $2.3490\cdot10^{-5}$ is achieved at epoch 1284, {and the corresponding test error is  $3.7580\cdot10^{-5}$}.
This selected model achieves a test error of $3.7580\cdot10^{-5}$ on the test data set,
thus we achieve a better accuracy compared to the non-convex scenario: both the best validation loss and the test error are 1 order of magnitude smaller.
This finding is inline with the higher regularity that characterizes this dataset, with respect to the one used for the non-convex case.
The hyperparameters and error values for this experiment are summarized in Table~\ref{tab:results_all}.

\begin{figure}[!t]
    \begin{subfigure}[b]{\linewidth}
\begin{tikzpicture}
\begin{axis}[enlargelimits=true,
width=\linewidth,
height=.45\linewidth,
ylabel={constant value},
xlabel={\#Convex polygon data - Poincar\'{e} inequality},
grid=major, grid style={dashed,gray!30},
grid style=dashed,
xmin = 0,
xmax = 385,
ymin=.2,
ymax=0.4,
legend style={at={(rel axis cs:0.5,-0.25)}, legend columns=6,fill=none,draw=none,anchor=north, /tikz/every even column/.append style={column sep=12pt}},
]

\addplot+[cnrblue, every mark/.append style = {fill = cnrblue}, mark = *, mark size = 1.5pt, mark options = {solid}, only marks] table [x=, y=y, col sep=comma]{hypertuning/clean_train_convex_poincare_labels.csv};\addlegendentry{data}

\addplot+[color=red, thick, mark size = 0pt] coordinates {(-50,0.254980) (500,0.254980)};\addlegendentry{mean}

\addplot+[color=dgreen, thick, densely dashed, mark size = 0pt] coordinates {(-50,0.26857) (500,0.26857)};\addlegendentry{$+$std}

\addplot+[color=olive, thick, densely dashed, mark size = 0pt] coordinates {(-50,0.24139) (500,0.24139)};\addlegendentry{$-$std}

\addplot+[color=cyan, thick, densely dashed, mark size = 0pt] coordinates {(-50,0.221704) (500,0.221704)};\addlegendentry{min}

\addplot+[color=magenta, thick, densely dashed, mark size = 0pt] coordinates {(-50,0.296813) (500,0.296813)};\addlegendentry{max}

\end{axis}
\end{tikzpicture}\caption{}
    \end{subfigure}
    \begin{subfigure}[b]{\linewidth}
\begin{tikzpicture}
\begin{axis}[enlargelimits=true, ymode = log,
width=\linewidth,
height=.3\linewidth,
xlabel={epochs},
ylabel={MSE},
grid=major, grid style={dashed,gray!30},
grid style=dashed,
xmin=0, xmax=1500,
ymax=0.1,
ymin=0.00001,
ytick={1e-1, 1e-2, 1e-3, 1e-4, 1e-5},
yticklabels={$10^{-1}$, $10^{-2}$, $10^{-3}$, $10^{-4}$, $10^{-5}$}]

\addplot+[cnrblue, every mark/.append style = {fill = cnrblue}, mark = *, mark size = 0.5pt, mark options = {solid}, opacity = 0.65] table [x=, y=trainloss, col sep=comma]{hypertuning/convex/poincare/hypertuning_convex_poincare_L_1_N_64_lr_0.005_losses.csv};\addlegendentry{training}

\addplot+[orange, every mark/.append style = {fill = orange}, mark size = 0.5pt, mark options = {solid}, opacity = 0.35] table [x=, y=valloss, col sep=comma]{hypertuning/convex/poincare/hypertuning_convex_poincare_L_1_N_64_lr_0.005_losses.csv};\addlegendentry{validation}
\end{axis}
\end{tikzpicture}\caption{}
    \end{subfigure}
    \caption{(a) Statistical analysis of the training target labels for learning the Poincar\'{e} constant on convex polygons: the data (blue), the mean value (read), the standard deviation (green and olive), the minimum (cyan) and the maximum values (magenta).
    (b) Training and validation loss over epochs for Poincar\'{e} inequality on convex polygons.
    }\label{fig:convex_poincare}
\end{figure}

\subsection{Inverse inequality on non-convex polygons.}
The training data for the inverse inequality on non-convex polygons
exhibit a mean of $10.5392$ with a standard deviation of  $4.2319$, indicating a moderate level of variability around the central value. Observed values range from $5.4543$ to $28.9294$, yielding a total range of $23.4751$. This broader spread reflects greater dispersion around the mean, suggesting the presence of more diverse or extreme observations. The data statistics are illustrated in Figure~\ref{fig:nonconvex_inverse}~(b).
The relatively large standard deviation -- about 40{\%} of the mean -- highlights a less homogeneous distribution compared to the Poincar\'{e} inequality case.

The distribution of the data affects also the model selection and its performance. More precisely, 
The complexity of the training set leds to a MLP architecture, characterized by $7$ hidden layers and $1024$ units per layer, to be trained via Adam optimizer with learning rate $10^{-2}$. { This architecture is the bigger one selected, among all the considered inequalities and polygon configurations, see for comparison Table~\ref{tab:results_all}.}
Moreover, the best validation loss value of $4.104611\cdot10^{-1}$ is achieved after 1139 epochs and this model leads to a test error of $5.193515\cdot10^{-1}$. The training and validation loss are illustrated in Figure~\ref{fig:nonconvex_inverse}~(b); { The model does not tend to overfit, though both training and validation error stagnates in the range $[0.1,1]$.}

\begin{figure}[!t]
\begin{subfigure}[b]{\linewidth}
\begin{tikzpicture}
\begin{axis}[enlargelimits=true,
width=\linewidth,
height=.45\linewidth,
ylabel={constant value},
xlabel={\#Non-convex polygon data - Inverse inequality},
grid=major, grid style={dashed,gray!30},
grid style=dashed,
xmin = 0,
xmax = 385,
ymin=5,
ymax=30,
legend style={at={(rel axis cs:0.5,-0.25)}, legend columns=6,fill=none,draw=none,anchor=north, /tikz/every even column/.append style={column sep=12pt}},
]

\addplot+[cnrblue, every mark/.append style = {fill = cnrblue}, mark = *, mark size = 1.5pt, mark options = {solid}, only marks] table [x=, y=y, col sep=comma]{hypertuning/clean_train_nonconvex_inverse_labels.csv};\addlegendentry{data}

\addplot+[color=red, thick, mark size = 0pt] coordinates {(-50,10.539236) (500,10.539236)};\addlegendentry{mean}

\addplot+[color=dgreen, thick, densely dashed, mark size = 0pt] coordinates {(-50,14.7710940) (500,14.7710940)};\addlegendentry{$+$std}

\addplot+[color=olive, thick, densely dashed, mark size = 0pt] coordinates {(-50,6.3073780) (500,6.3073780)};\addlegendentry{$-$std}

\addplot+[color=cyan, thick, densely dashed, mark size = 0pt] coordinates {(-50,5.454347) (500,5.454347)};\addlegendentry{min}

\addplot+[color=magenta, thick, densely dashed, mark size = 0pt] coordinates {(-50,28.929404) (500,28.929404)};\addlegendentry{max}

\end{axis}
\end{tikzpicture}\caption{}
\end{subfigure}
\begin{subfigure}[b]{\linewidth}
\begin{tikzpicture}
\begin{axis}[enlargelimits=true, ymode = log,
width=\linewidth,
height=.3\linewidth,
xlabel={epochs},
ylabel={MSE},
grid=major, grid style={dashed,gray!30},
grid style=dashed,
xmin=0, xmax=1500,
ymax=100,
ymin=0.01,
ytick={1e2, 1e1, 1, 1e-1, 1e-2},
yticklabels={$10^{2}$, $10$, $1$, $10^{-1}$, $10^{-2}$}]

\addplot+[cnrblue, every mark/.append style = {fill = cnrblue}, mark = *, mark size = 0.5pt, mark options = {solid}, opacity = 0.65] table [x=, y=trainloss, col sep=comma]{hypertuning/nonconvex/inverse/hypertuning_nonconvex_inverse_L_7_N_1024_lr_0.01_losses.csv};\addlegendentry{training}

\addplot+[orange, every mark/.append style = {fill = orange}, mark size = 0.5pt, mark options = {solid}, opacity = 0.35] table [x=, y=valloss, col sep=comma]{hypertuning/nonconvex/inverse/hypertuning_nonconvex_inverse_L_7_N_1024_lr_0.01_losses.csv};\addlegendentry{validation}
\end{axis}
\end{tikzpicture}\caption{}
\end{subfigure}
    \caption{(a) Statistical analysis of the training target labels for learning the Inverse constant on non-convex polygons: the data (blue), the mean value (red), the standard deviation (green and olive), the minimum (cyan) and the maximum values (magenta).
    (b) Training and validation loss over 1500 epochs for Inverse inequality on non-convex polygons.}\label{fig:nonconvex_inverse}
\end{figure}

\subsection{Inverse inequality on convex polygons}
A similar trend appears on convex polygons, where the inverse inequality data has a mean of $7.8873$ and a standard deviation of $2.3920$ -- over 30{\%} of the mean. Values span from $4.6239$ to $19.2646$, indicating a smaller spread than the non-convex data. This suggests decreased variability and a potentially more homogeneous distribution of values;
the data statistics are illustrated in Figure~\ref{fig:convex_inverse}~(a).
The hyperparameter optimization leads to an MLP, { smaller than the one selected in the non-convex scenario}, with 3 hidden layers and 256 units per layer, to be trained via Adam optimizer with learning rate $10^{-3}$; the training and validation loss are illustrated in Figure~\ref{fig:convex_inverse}~(b).
In the end, the best validation loss value of $3.0341\cdot10^{-2}$ is achieved after 548 epochs and this model leads to a test error of $9.1029\cdot10^{-1}$. { Differently from the other two inequalities, we register a lower test error for the inverse inequality on non-convex polygons intead of on convex ones, see Table~\ref{tab:results_all} for further comparison.}

\begin{figure}[!t]
\begin{subfigure}[b]{\linewidth}
\begin{tikzpicture}
\begin{axis}[enlargelimits=true,
width=\linewidth,
height=.45\linewidth,
ylabel={constant value},
xlabel={\#Convex polygon data - Inverse inequality},
grid=major, grid style={dashed,gray!30},
grid style=dashed,
xmin = 0,
xmax = 385,
ymin=5,
ymax=30,
legend style={at={(rel axis cs:0.5,-0.25)}, legend columns=6,fill=none,draw=none,anchor=north, /tikz/every even column/.append style={column sep=12pt}},
]

\addplot+[cnrblue, every mark/.append style = {fill = cnrblue}, mark = *, mark size = 1.5pt, mark options = {solid}, only marks] table [x=, y=y, col sep=comma]{hypertuning/clean_train_convex_inverse_labels.csv};\addlegendentry{data}

\addplot+[color=red, thick, mark size = 0pt] coordinates {(-50,7.887386) (500,7.887386)};\addlegendentry{mean}

\addplot+[color=dgreen, thick, densely dashed, mark size = 0pt] coordinates {(-50,10.2793379) (500,10.2793379)};\addlegendentry{$+$std}

\addplot+[color=olive, thick, densely dashed, mark size = 0pt] coordinates {(-50,5.4954340) (500,5.495434)};\addlegendentry{$-$std}

\addplot+[color=cyan, thick, densely dashed, mark size = 0pt] coordinates {(-50,4.623949) (500,4.623949)};\addlegendentry{min}

\addplot+[color=magenta, thick, densely dashed, mark size = 0pt] coordinates {(-50,19.2646143) (500,19.264614)};\addlegendentry{max}

\end{axis}
\end{tikzpicture}\caption{}
\end{subfigure}

\begin{subfigure}[b]{\linewidth}
\begin{tikzpicture}
\begin{axis}[enlargelimits=true, ymode = log,
width=\linewidth,
height=.3\linewidth,
xlabel={epochs},
ylabel={MSE},
grid=major, grid style={dashed,gray!30},
grid style=dashed,
xmin=0, xmax=1500,
ymax=100,
ymin=0.01,
ytick={1e2, 1e1, 1, 1e-1, 1e-2},
yticklabels={$10^{2}$, $10$, $1$, $10^{-1}$, $10^{-2}$}]

\addplot+[cnrblue, every mark/.append style = {fill = cnrblue}, mark = *, mark size = 0.5pt, mark options = {solid}, opacity = 0.65] table [x=, y=trainloss, col sep=comma]{hypertuning/convex/inverse/hypertuning_convex_inverse_L_3_N_256_lr_0.001_losses.csv};\addlegendentry{training}

\addplot+[orange, every mark/.append style = {fill = orange}, mark size = 0.5pt, mark options = {solid}, opacity = 0.35] table [x=, y=valloss, col sep=comma]{hypertuning/convex/inverse/hypertuning_convex_inverse_L_3_N_256_lr_0.001_losses.csv};\addlegendentry{validation}
\end{axis}
\end{tikzpicture}\caption{}
\end{subfigure}

\caption{(a) Statistical analysis of the training target labels for learning the Inverse constant on convex polygons: the data (blue), the mean value (red), the standard deviation (green and olive), the minimum (cyan) and the maximum values (magenta).
    (b) Training and validation loss over 1500 epochs for Inverse inequality on convex polygons.}
\label{fig:convex_inverse}
\end{figure}

\subsection{Trace inequality on non-convex polygons}
As concerns the trace inequality, in the non-convex case, the training data shows a mean of $5.2002$ and a standard deviation of $1.1705$. The observed values span from $3.2432$ to $12.3412$, yielding a total range of $9.0980$. While the variability remains moderate, the broad spread suggests a more diverse distribution -- possibly due to the geometric complexity of non-convex domains. The training labels are illustrated in Figure~\ref{fig:nonconvex_trace}~(a).
{To predict the geometric constant, we selected a MLP of moderate size, characterized by $3$ hidden layers and $32$ units per layer. The training is then performed via Adam optimized with learning rate $5\cdot10^{-3}$. The training and validation loss are shown in Figure~\ref{fig:nonconvex_trace}~(b). More precisely the best validation loss of $5.7048\cdot10^{-2}$ is achieved after 1209 epochs and the corresponding model registers a test error of $5.5408\cdot10^{-2}$.}

\begin{figure}[!t]
    \begin{subfigure}[b]{\linewidth}
\begin{tikzpicture}
\begin{axis}[enlargelimits=true,
width=\linewidth,
height=.45\linewidth,
ylabel={constant value},
xlabel={\#Non-convex polygon data - Trace inequality},
grid=major, grid style={dashed,gray!30},
grid style=dashed,
xmin = 0,
xmax = 385,
ymin=3,
ymax=10,
legend style={at={(rel axis cs:0.5,-0.25)}, legend columns=6,fill=none,draw=none,anchor=north, /tikz/every even column/.append style={column sep=12pt}},
]

\addplot+[cnrblue, every mark/.append style = {fill = cnrblue}, mark = *, mark size = 1.5pt, mark options = {solid}, only marks] table [x=, y=y, col sep=comma]{hypertuning/clean_train_nonconvex_trace_labels.csv};\addlegendentry{data}

\addplot+[color=red, thick, mark size = 0pt] coordinates {(-50,5.258960) (500,5.258960)};\addlegendentry{mean}

\addplot+[color=dgreen, thick, densely dashed, mark size = 0pt] coordinates {(-50,6.2752970) (500,6.2752970)};\addlegendentry{$+$std}

\addplot+[color=olive, thick, densely dashed, mark size = 0pt] coordinates {(-50,4.2426230) (500,4.2426230)};\addlegendentry{$-$std}

\addplot+[color=cyan, thick, densely dashed, mark size = 0pt] coordinates {(-50,3.663829) (500,3.663829)};\addlegendentry{min}

\addplot+[color=magenta, thick, densely dashed, mark size = 0pt] coordinates {(-50,9.668667) (500,9.668667)};\addlegendentry{max}

\end{axis}
\end{tikzpicture}\caption{}
    \end{subfigure}

 \begin{subfigure}[b]{\linewidth}
\begin{tikzpicture}
\begin{axis}[enlargelimits=true, ymode = log,
width=\linewidth,
height=.3\linewidth,
xlabel={epochs},
ylabel={MSE},
grid=major, grid style={dashed,gray!30},
grid style=dashed,
xmin=0, xmax=1500,
ymax=20,
ymin=0.01,
ytick={1e1, 1, 1e-1, 1e-2},
yticklabels={$10$, $1$, $10^{-1}$, $10^{-2}$}]

\addplot+[cnrblue, every mark/.append style = {fill = cnrblue}, mark = *, mark size = 0.5pt, mark options = {solid}, opacity = 0.65] table [x=, y=trainloss, col sep=comma]{hypertuning/nonconvex/trace/hypertuning_nonconvex_trace_L_3_N_32_lr_0.005_losses.csv};\addlegendentry{training}

\addplot+[orange, every mark/.append style = {fill = orange}, mark size = 0.5pt, mark options = {solid}, opacity = 0.35] table [x=, y=valloss, col sep=comma]{hypertuning/nonconvex/trace/hypertuning_nonconvex_trace_L_3_N_32_lr_0.005_losses.csv};\addlegendentry{validation}
\end{axis}
\end{tikzpicture}\caption{}
    \end{subfigure}
    
    \caption{(a) Statistical analysis of the training target labels for learning the Trace constant on non-convex polygons: the data (blue), the mean value (red), the standard deviation (green and olive), the minimum (cyan) and the maximum values (magenta).
    (b) Training and validation loss over 1500 epochs for Trace inequality on non-convex polygons.
    }\label{fig:nonconvex_trace}
\end{figure}

\subsection{Trace inequality on convex polygons.}
Finally, the training data for the trace inequality on convex polygons has a mean of $3.8074$ and a standard deviation of $0.3679$, indicating relatively low variability around the average. The data range from a minimum of $3.3071$ to a maximum of $5.5320$, giving a total range of $2.2249$. With the standard deviation amounting to less than 10{\%}
of the mean, the distribution is fairly concentrated, suggesting a consistent pattern with limited deviation from the central tendency.
{The optimized architecture has $4$ hidden layers more than the one of the convex case, with the same number of hidden units per layer, namely $32$; also the learning rate employed for the training procedure is the same: $5\cdot10^{-3}$. The best validation loss of $3.1301\cdot10^{-2}$ is achieved at epoch 1352 and the corresponding model has a test error of $7.5794\cdot10^{-2}$.
}

\begin{figure}[!t]
\begin{subfigure}[b]{\linewidth}
\begin{tikzpicture}
\begin{axis}[enlargelimits=true,
width=\linewidth,
height=.45\linewidth,
ylabel={constant value},
xlabel={\#Convex polygon data - Trace inequality},
grid=major, grid style={dashed,gray!30},
grid style=dashed,
xmin = 0,
xmax = 385,
ymin=3,
ymax=10,
legend style={at={(rel axis cs:0.5,-0.25)}, legend columns=6,fill=none,draw=none,anchor=north, /tikz/every even column/.append style={column sep=12pt}},
]

\addplot+[cnrblue, every mark/.append style = {fill = cnrblue}, mark = *, mark size = 1.5pt, mark options = {solid}, only marks] table [x=, y=y, col sep=comma]{hypertuning/clean_train_convex_trace_labels.csv};\addlegendentry{data}

\addplot+[color=red, thick, mark size = 0pt] coordinates {(-50,4.177908) (500,4.177908)};\addlegendentry{mean}

\addplot+[color=dgreen, thick, densely dashed, mark size = 0pt] coordinates {(-50,4.775170) (500,4.775170)};\addlegendentry{$+$std}

\addplot+[color=olive, thick, densely dashed, mark size = 0pt] coordinates {(-50,3.5806460) (500,3.5806460)};\addlegendentry{$-$std}

\addplot+[color=cyan, thick, densely dashed, mark size = 0pt] coordinates {(-50,3.141536) (500,3.141536)};\addlegendentry{min}

\addplot+[color=magenta, thick, densely dashed, mark size = 0pt] coordinates {(-50,6.909680) (500,6.909680)};\addlegendentry{max}

\end{axis}
\end{tikzpicture}
\end{subfigure}

\begin{subfigure}[b]{\linewidth}
\begin{tikzpicture}
\begin{axis}[enlargelimits=true, ymode = log,
width=\linewidth,
height=.3\linewidth,
xlabel={epochs},
ylabel={MSE},
grid=major, grid style={dashed,gray!30},
grid style=dashed,
xmin=0, xmax=1500,
ymax=20,
ymin=0.01,
ytick={1e1, 1, 1e-1, 1e-2},
yticklabels={$10$, $1$, $10^{-1}$, $10^{-2}$}]

\addplot+[cnrblue, every mark/.append style = {fill = cnrblue}, mark = *, mark size = 0.5pt, mark options = {solid}, opacity = 0.65] table [x=, y=trainloss, col sep=comma]{hypertuning/convex/trace/hypertuning_convex_trace_L_7_N_32_lr_0.005_losses.csv};\addlegendentry{training}

\addplot+[orange, every mark/.append style = {fill = orange}, mark size = 0.5pt, mark options = {solid}, opacity = 0.35] table [x=, y=valloss, col sep=comma]{hypertuning/convex/trace/hypertuning_convex_trace_L_7_N_32_lr_0.005_losses.csv};\addlegendentry{validation}
\end{axis}
\end{tikzpicture}
\end{subfigure}

\caption{(a) Statistical analysis of the training target labels for learning the Trace constant on convex polygons: the data (blue), the mean value (red), the standard deviation (green and olive), the minimum (cyan) and the maximum values (magenta).
    (b) Training and validation loss over 1500 epochs for Trace inequality on convex polygons.}\label{fig:convex_trace}
\end{figure}

\subsection{Dealing with arbitrary polygon sizes}    
 \label{general_polygons} 
{In this section, we illustrate how the proposed method can be applied not only to polygons of unitary diameter, but also to more general ones, characterized by an arbitrary diameter value.}
Let us consider for example, the Poincar\'{e} inequality on a domain of arbitrary size. It is well known that this inequality takes the form:
\begin{equation}
\label{Poincare_h}
	\norm{w}_{0,\om}\leq \Cp \, h \norm{\nabla w}_{0, \om} \qquad \forall \ w \in V,
\end{equation}
where $h = \mathrm{diam}(\om)$ denotes the diameter of $\om$, and $C_P$ is a constant independent of $h$.

In order to apply the method to polygons of any size,
for instance element of a mesh as in Figure \ref{fig:polygons}, for each polygon $K$, 
we scale the related geometric constant with the following factor
 \[
     \delta = \delta\left(K\right) := \frac{1}{\mathrm{diam}\left(K\right)}.
 \]
 Therefore, the rescaled input attributes are determined as follows:
 \begin{align*}
 \widetilde{\CC}  & = \delta \CC  & \widetilde{\IC}  & = \delta \IC  & \widetilde{\CR}  & = \CR          & \widetilde{\AR} & = \delta^2 \AR\\
 \widetilde{\KE}  & = \delta^2\KE & \widetilde{\KAR} & = \KAR        & \widetilde{\APR} & = \APR         & \widetilde{\SE} & = \delta\SE\\
 \widetilde{\sSE} & = \sSE        & \widetilde{\ER}  & = \ER         & \widetilde{\MPD} & = \delta \MPD  & \widetilde{\sMPD} & = \sMPD\\ 
 \widetilde{\MA}  & = \MA         & \widetilde{\MX}  & = \MX         & \widetilde{\SR}  & = \SR          & \widetilde{\ISO} & = \ISO.
 \end{align*}
 Note that the linear variables are scaled by a factor $\delta$, that cancels out for the variables defined as the ratio of two linear ones, \eg $\widetilde{\CR}$; the areas are scaled by a square factor $\delta^2$; the angles are not affected by the length deformation; the isotropy becomes the ratio between minimum and maximum eigenvalue of the covariance matrix $\widetilde{M}_{cov}\coloneq\frac{1}{\widetilde{\AR}}\int_{\widetilde{K}}(x-\bar{x}_{\widetilde{K}})^T(x-\bar{x}_{\widetilde{K}})dx$, where 
 $\widetilde{K} = \delta K$,
 and with a change of variables we can prove that $\widetilde{M}_{cov}=M_{cov}$.

The considered attributes are given as input to our ANN tool:
as before each neuron in the input layer corresponds to a different geometrical feature of a polygon $K$ and 
the output layer consists of a single linear unit, providing the estimate of the corresponding considered \emph{scaled} constant.
To obtain the constant corresponding to $K$, we need to rescale the output by a factor $\frac{1}{\delta_\om}$.
The summary of the proposed method is illustrated in Figure~\ref{rescaledANN}.

	\begin{figure}[!t]
		\centering

            \begin{tikzpicture}[
    node distance=1.5cm and 1.5cm,
    every node/.style={align=center},
    arrow/.style={black, -{Latex[round]}, thick},
    rescale/.style={black, font=\scriptsize},
    neuron/.style={circle, draw, minimum size=6pt, inner sep=0pt},
    layer/.style={matrix of nodes, nodes=neuron, column sep=4mm, row sep=4mm, anchor=center},
    stage/.style={draw=cnrblue!50, fill=cnrblue!10, thick, rounded corners, inner sep=0.5em}
  ]
  
  \node (X) {$K$};
  
  \node[right=of X] (Xt) {$\mathbf{m}$};

  \draw[arrow] (X) -- (Xt);
  
  \node (scale1) [rescale, below=0.1cm of $(X)!0.5!(Xt)$] {Rescaling ($\delta$)};

  \node (output-1-1) [right = of Xt] {
    \begin{tikzpicture}[scale=0.5]
					\tikzstyle{unit}=[draw, shape=circle, minimum size=0.035cm]
					
					\tikzstyle{hidden}=[draw,shape=circle, minimum size=0.035cm]
					
					\node[unit](x0) at (0,3.5){};
					\node at (0,2.8){$\vdots$};
					\node[unit](x1) at (0, 1.75){};
					\node at (0,1){$\vdots$};
					\node[unit](xd) at (0,0){};
					
					\node[hidden](h10) at (3-1.5,4){};
					\node[hidden](h11) at (3-1.5,2.5){};
					\node at (3-1.5,1.5){$\vdots$};
					\node[hidden](h1m) at (3-1.5,-0.5){};
					
					\node(h22) at (5-1.5,0){};
					\node(h21) at (5-1.5,2){};
					\node(h20) at (5-1.5,4){};
					
					\node(d3) at (6-2,0){$\dots$};
					\node(d2) at (6-2,2){$\dots$};
					\node(d1) at (6-2,4){$\dots$};
					
					\node(hL12) at (7-1.5,0){};
					\node(hL11) at (7-1.5,2){};
					\node(hL10) at (7-1.5,4){};
					
					\node[hidden](hL0) at (9-1.5,4){};
					\node[hidden](hL1) at (9-1.5,2.5){};
					\node at (9-1.5,1.5){$\vdots$};
					\node[hidden](hLm) at (9-1.5,-0.5){};

					\node[unit](y2) at (11-2,2){};

					\draw[-latex] (x0) -- (h11);
					\draw[-latex] (x0) -- (h1m);
					
					\draw[-latex] (x1) -- (h11);
					\draw[-latex] (x1) -- (h1m);
					
					\draw[-latex] (xd) -- (h11);
					\draw[-latex] (xd) -- (h1m);

					\draw[-latex] (hL0) -- (y2);

					\draw[-latex] (hL1) -- (y2);

					\draw[-latex] (hLm) -- (y2);

					\draw[-latex, path fading=east] (h10) -- (h21) node[pos=0.3, sloped, above]  {};
					\draw[-latex, path fading=east] (h10) -- (h22);
					
					\draw[-latex, path fading=east] (h11) -- (h21);
					\draw[-latex, path fading=east] (h11) -- (h22);
					
					\draw[-latex, path fading=east] (h1m) -- (h21);
					\draw[-latex, path fading=east] (h1m) -- (h22);
					
					\draw[-latex, path fading=west] (hL10) -- (hL1);
					\draw[-latex, path fading=west] (hL11) -- (hL1);
					\draw[-latex, path fading=west] (hL12) -- (hL1);
					
					\draw[-latex, path fading=west] (hL10) -- (hLm);
					\draw[-latex, path fading=west] (hL11) -- (hLm);
					\draw[-latex, path fading=west] (hL12) -- (hLm);

				\end{tikzpicture}
    
    };

  \node[right=1.2cm of output-1-1] (yt) {$C(\mathbf{m})$};
  \node[right=of yt] (yp) {$C(K)$};

  \draw[arrow] (Xt) -- (output-1-1);
  \draw[arrow] (output-1-1) -- (yt);
  \draw[arrow] (yt) -- (yp);
  \node (scale2) [rescale, below=0.1cm of $(yt)!0.5!(yp)$] {Rescaling ($\nicefrac{1}{\delta}$)};
  
\begin{scope}[on background layer]
  \node[stage, fit=(X)(Xt)(scale1)] (inputbox) {};
  \node[above=0.2cm of inputbox] {Input Features};

  \node[stage, fit=(output-1-1)] (nnbox) {};
  \node[above=0.2cm of nnbox] {Neural Network};

  \node[stage, fit=(yt)(yp)(scale2)] (outbox) {};
  \node[above=0.2cm of outbox] {Estimate}; 
\end{scope}

\end{tikzpicture}
        
		\caption{The geometric features are computed for a polygon $K$, they are rescaled by a factor $\delta$ and passed to the ANN model for training. In the evaluation phase only $\mathbf{m}$ is passed to the ANN and the output $C(\mathbf{m})$ is rescaled by a factor $\frac{1}{\delta}$ to obtain the evaluation of the considered constant corresponding to $K$.}
		\label{rescaledANN}
	\end{figure}

\section{Conclusions and further developments}\label{sec:6}
{We propose a method for estimating three constants appearing in key inequalities on generic polygons, based on Artificial Neural Networks (ANNs) that take as input a set of geometric features. The loss function on the validation set typically converges within a few thousand iterations. For each constant, the trained model achieves an accuracy of approximately two significant digits, representing a substantial improvement over the overly pessimistic upper bounds commonly used in practice.

The optimal ANN architecture is relatively shallow, and the evaluation of the constants after training—which is performed offline once and for all—involves only simple operations such as additions and multiplications. As a result, the online phase of the method is extremely fast.
Additionally, the required dataset for training is small, significantly reducing not only the offline training time but also the effort needed for data generation and preprocessing.

Notably, the same procedure can be applied across the different constants with minimal modifications, demonstrating the flexibility of the approach. This indicates  that the ANN-based framework can be readily extended to estimate other geometric constants beyond those considered.

A particularly interesting extension concerns the three-dimensional setting. Our metric-based representation appears especially effective in 3D, as it characterizes a polytope without requiring its explicit geometric definition—which is often computationally expensive—and instead relies on fixed, measurable quantities intrinsic to the shape.

Future work will explore the extension of this method to the computation of constants for polytopes in three dimensions, as well as for polygons arising from general polygonal meshes, with the aim of supporting the design and analysis of polytopal methods for partial differential equations (PDEs).
}

\section*{Acknowledgments}
This paper has been co–funded by the MUR Progetti di Ricerca di Rilevante Interesse Nazionale (PRIN) Bando 2022/PNRR (grant P2022BH5CB, NextGenerationEU). 
SI and MP have been partially supported by ICSC-Centro Nazionale di Ricerca in High Performance Computing, Big Data, and Quantum Computing funded by the European Union–NextGenerationEU plan.
The authors are members of INDAM-GNCS.

\appendix
\section{Non-convex polygon generation}

\begin{figure}[!ht]
		\centering\includegraphics[width=.7\linewidth]{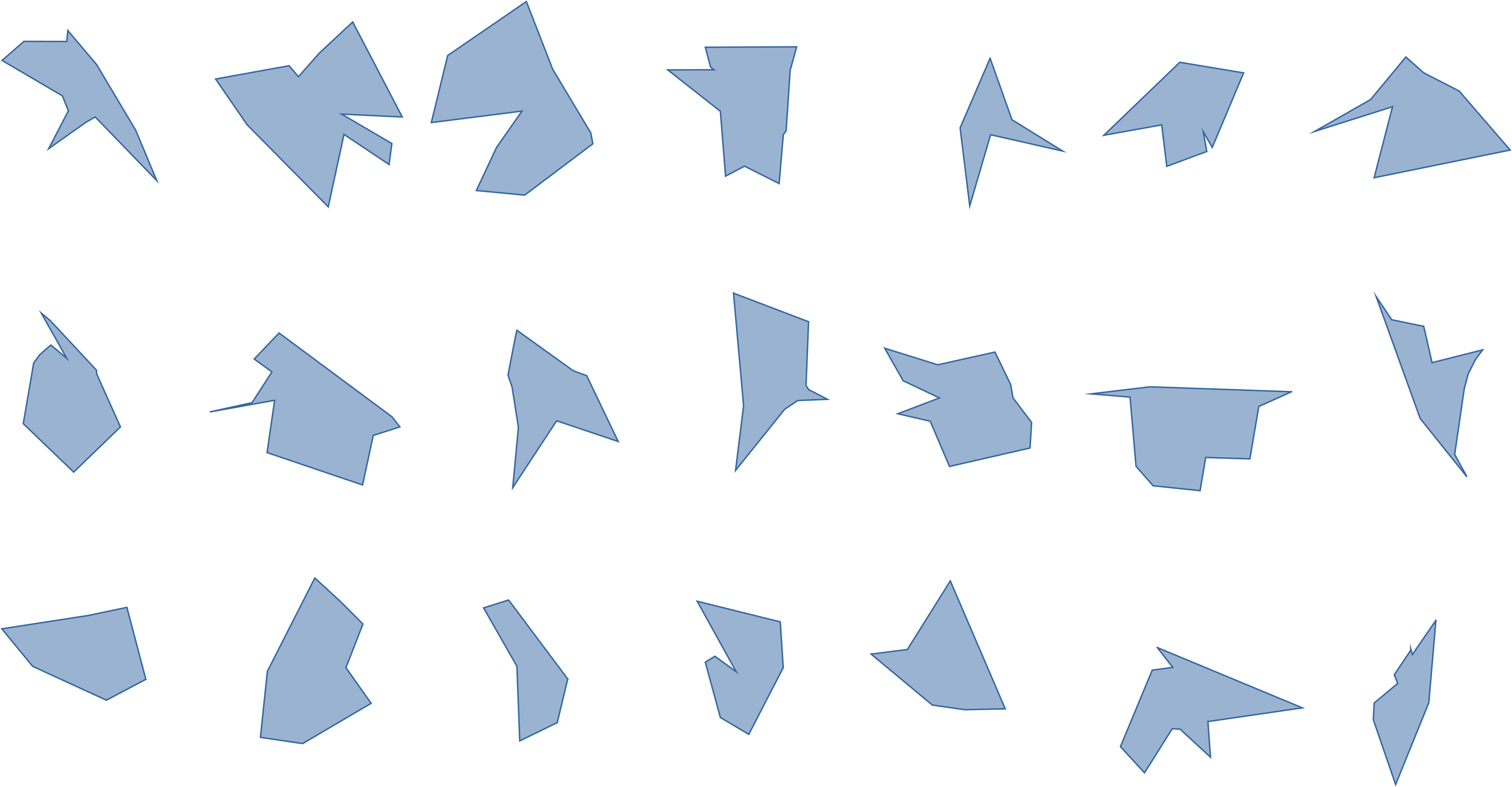}
		\caption{Example of 21 non-convex polygonal domains, generated with Algorithm~\ref{alg:nonconvex-gen}.}
		\label{fig:non-convex}
	\end{figure}

We describe the procedure considered 
to generate random non-convex polygons. 
For any polygon, we generate random angle values sampled in $[\nicefrac{\pi}{6}, \pi]$,
such that their sum is equal to $2\pi$.
For each random angle $\theta$, we associate a random length $\ell_\theta$. Subsequently, we define the vertices of the polygon in terms of polar coordinates with center in $(0,0)$ and parameters $\ell_\theta$ and cumulative sum of $\theta$. The pseudocode is in Algorithm\ref{alg:nonconvex-gen}, whereas an illustrative example of 21 random non-convex polygons is shown in Figure~\ref{fig:non-convex}.
\begin{algorithm2e}[!ht]
	\SetAlgoLined
	\KwIn{Number of vertices $n\in\mathbb{N}$}
        \For{$i=1, \ldots, n$}{
        Sample the angles $\theta_i\sim\mathcal{U}[\nicefrac{\pi}{6}, \pi]$\\
        Sample random lengths $\ell_i\sim\mathcal{U}[0,1]$\\
        }
        Set $\Theta = \sum_{i=1}^n\theta_i$\\
        \For{$i=1, \ldots, n$}{
        Normalize the angles in $[0,2\pi]$:\\
        $\theta_i = \nicefrac{\theta_i}{\Theta}2\pi$\\
        }
        Set $\alpha = 0$\\
        \For{$i=1, \ldots, n$}{
        $\alpha = \alpha + \theta_i$\\
        $x_i = \ell_i\cos(\alpha)$\\
        $y_i = \ell_i\sin(\alpha)$\\
        Append $(x_i, y_i)$ to $V$
        }
    \KwOut{The set of vertices $V$.}
    \caption{Non-convex polygons generation.}
    \label{alg:nonconvex-gen}
\end{algorithm2e}
    
\bibliographystyle{elsarticle-num-names}
\bibliography{bibliography}

@article{eigenvalues,
	author = {Matculevich, S. and Repin, S.},
	title = {Sharp bounds of constants in Poincare type inequalities for polygonal domains},
	journal = {arXiv},
	year ={ 2015}
}

@article{attene2021,
	title = {Benchmarking the geometrical robustness of a Virtual Element Poisson solver},
	journal = {Mathematics and Computers in Simulation},
	volume = {190},
	pages = {1392-1414},
	year = {2021},
	author = {Attene, Marco and Biasotti, Silvia and Bertoluzza, Silvia and Cabiddu, Daniela and Livesu, Marco and Patanè, Giuseppe and Pennacchio, Micol and Prada, Daniele and Spagnuolo, Michela}
}

@misc{WinNT,
    author={Wolf Vollprecht},
	  title = {Scikit-geometry Documentation},
    howpublished = {\url{https://scikit-geometry.github.io/scikit-geometry/}},
	  note = {Accessed: 2025-10-14},
    year={2019},
}

@INPROCEEDINGS{tensorflow,
	  author={Ertam, Fatih and Aydın, Galip},
	  booktitle={2017 International Conference on Computer Science and Engineering (UBMK)}, 
	  title={Data classification with deep learning using {T}ensorflow}, 
	  year={2017},
	  pages={755-758},
	  doi={10.1109/UBMK.2017.8093521}
}

@article{isosceles_cp,
	title = {Minimizing {N}eumann fundamental tones of triangles: An optimal {P}oincaré inequality},
	journal = {Journal of Differential Equations},
	volume = {249},
	number = {1},
	pages = {118-135},
	year = {2010},
	author = {R.S. Laugesen and B.A. Siudeja}
}

@article{exact_cp,
	author = {Nazarov, A. I. and Repin, S. I.},
	title = {Exact constants in {P}oincaré type inequalities for functions with zero mean boundary traces},
	journal = {Mathematical Methods in the Applied Sciences},
	volume = {38},
	number = {15},
	pages = {3195-3207},
	year = {2015}
}

@article{convex_cp,
	author={Payne, L. E. and Weinberger, H. F.},
	year={1960},
	title={An optimal {P}oincaré inequality for convex domains},
	journal={Archive for Rational Mechanics and Analysis},
	pages={286-292},
	volume={5}
	
}

@article{repin,
	title={Computable majorants of constants in the {P}oincar\'e and {F}riedrichs inequalities},
	author={Repin, S.},
	journal={Journal of Mathematical Sciences},
	volume={186},
	number={2},
	year={2012}
}

@book{lame,
	title={Le{\c{c}}ons sur la th{\'e}orie math{\'e}matique de l'{\'e}lasticit{\'e} des corps solides},
	author={Lam{\'e}, Gabriel},
	year={1852},
	publisher={Bachelier}
}

@article{equilateral_cp,
	title={Eigenstructure of the equilateral triangle, Part II: The Neumann problem},
	author={McCartin, Brian J},
	journal={Mathematical Problems in Engineering},
	volume={8},
	number={6},
	pages={517--539},
	year={2002},
	publisher={Hindawi}
}

@book{brenner_scott,
  author    = {Brenner, Susanne C. and Scott, Ridgway},
  title     = {The Mathematical Theory of Finite Element Methods},
  publisher = {Springer},
  edition   = {3rd},
  year      = {2008}
}

@book{ern_guermond,
  author    = {Ern, Alexandre and Guermond, Jean-Luc},
  title     = {Theory and Practice of Finite Elements},
  publisher = {Springer},
  year      = {2004}
}

@article{ainsworth_inverse,
  author    = {Ainsworth, Mark and Coyle, J.},
  title     = {Hierarchic finite element bases on unstructured tetrahedral meshes},
  journal   = {International Journal for Numerical Methods in Engineering},
  volume    = {58},
  number    = {14},
  pages     = {2103--2130},
  year      = {2003}
}

@article{hughes_inequalities,
  author    = {Hughes, Thomas J.R. and Balestra, M. and Larson, M.G. and Masud, A.},
  title     = {What are C and h? Inequalities for the analysis and design of finite element methods},
  journal   = {Computer Methods in Applied Mechanics and Engineering},
  volume    = {197},
  number    = {49-50},
  pages     = {4855--4863},
  year      = {2008},
  doi       = {10.1016/j.cma.2008.06.008}
}

@book{toselli2005domain,
  title={Domain Decomposition Methods: Algorithms and Theory},
  author={Toselli, Andrea and Widlund, Olof},
  publisher={Springer},
  year={2005},
  series={Springer Series in Computational Mathematics},
  volume={34}
}

@book{steinbach2008numerical,
  title={Numerical Approximation Methods for Elliptic Boundary Value Problems: Finite and Boundary Elements},
  author={Steinbach, Olaf},
  publisher={Springer},
  year={2008}
}

@article{ainsworth2005posteriori,
  title={A posteriori error estimation for non-conforming finite element methods},
  author={Ainsworth, Mark},
  journal={SIAM Journal on Numerical Analysis},
  volume={42},
  number={6},
  pages={2320--2341},
  year={2005}
}

@book{verfurth2013posteriori,
  title={A Posteriori Error Estimation Techniques for Finite Element Methods},
  author={Verf{\"u}rth, R{\"u}diger},
  publisher={Oxford University Press},
  year={2013}
}

@article{beirao2013basic,
  title={Basic principles of virtual element methods},
  author={Beir{\~a}o da Veiga, L. and Brezzi, F. and Cangiani, A. and Manzini, G. and Marini, L.D. and Russo, A.},
  journal={Mathematical Models and Methods in Applied Sciences},
  volume={23},
  number={01},
  pages={199--214},
  year={2013}
}

@article{da2014mimetic,
  title={Mimetic finite difference methods for elliptic problems},
  author={Beir{\~a}o da Veiga, L. and Lipnikov, K. and Manzini, G.},
  journal={Foundations of Computational Mathematics},
  volume={14},
  number={6},
  pages={1147--1193},
  year={2014}
}

@article{sukumar2004finite,
  title={The natural element method in solid mechanics},
  author={Sukumar, N. and Tabarraei, A.},
  journal={International Journal for Numerical Methods in Engineering},
  volume={61},
  number={12},
  pages={2159--2181},
  year={2004}
}

@article{gain2014virtual,
  title={On the virtual element method for three-dimensional linear elasticity problems on arbitrary polyhedral meshes},
  author={Gain, A.L. and Talischi, C. and Paulino, G.H.},
  journal={Computer Methods in Applied Mechanics and Engineering},
  volume={282},
  pages={132--160},
  year={2014}
}

@article{yaro2017,
	title = {Error bounds for approximations with deep ReLU networks},
	journal = {Neural Networks},
	volume = {94},
	pages = {103-114},
	year = {2017},
	author = {Dmitry Yarotsky}
}

@book{goodfellow2016,
	title={Deep Learning},
	author={Ian Goodfellow and Yoshua Bengio and Aaron Courville},
	publisher={MIT Press},
	note={\url{http://www.deeplearningbook.org}},
	year={2016}
}

@book{aggarwal2018,
	title={Neural networks and deep learning},
	author={Aggarwal, Charu C.},
	year={2018},
	publisher={Springer}
}

\end{document}